\documentclass[12pt]{amsart}

\usepackage{amsfonts,amsmath,latexsym,amssymb,verbatim,amsbsy}
\usepackage{amsthm}

\usepackage{pstricks}

\theoremstyle{plain}
\newtheorem{THEOREM}{Theorem}[section]

\newtheorem{theorem}[THEOREM]{Theorem}
\newtheorem{corollary}[THEOREM]{Corollary}
\newtheorem{lemma}[THEOREM]{Lemma}

\theoremstyle{definition}

\newtheorem{definition}[THEOREM]{Definition}

\theoremstyle{remark}

\newcommand{\thm}[1]{Theorem~\ref{#1}}
\newcommand{\lem}[1]{Lemma~\ref{#1}}
\newcommand{\defin}[1]{Definition~\ref{#1}}

\newcommand{\N}{\ensuremath{\mathbb{N}}}   
\newcommand{\R}{\ensuremath{\mathbb{R}}}   
\newcommand{\C}{\ensuremath{\mathbb{C}}}   
\newcommand{\T}{\ensuremath{\mathbb{T}}}   

\def \a {\alpha}

\def \d {\delta}

\def \e {\varepsilon}
\def \f {\varphi}
\def \F {{\bf \Phi}}

\def \l {\lambda}

\def \n {\nabla}
\def \s {\sigma}
\def \th {\theta}
\def \Th {\Theta}

\def \fm {\mathfrak{m}}

\def \ba {{\bf a}}

\def \bB {{\bf B}}
\def \bC {{\bf C}}

\def \bP {{\bf P}}

\def \bI {{\bf I}}

\def \bE {{\bf E}}

\def \bG {{\bf G}}
\def \bT {{\bf T}}

\def \bX {{\bf X}}

\def \cL {\mathcal{L}}

\def \bFi {{\bf \Phi}}
\def \bPsi {{\bf \Psi}}

\def \M {\mathcal{M}}
\def \torus {{\T^n}}
\def \rn {\R^n \backslash\{0\}}

\def \lmax {\l_{\mathrm{max}}}
\def \lmin {\l_{\mathrm{min}}}

\def \Ms {Ma\~{n}e sequence}

\def \p {\partial}
\def \ra {\rightarrow}
\def \ss {\subset}

\newcommand{\der}[2]{(#1 \cdot \nabla) #2}

\DeclareMathOperator{\im}{Im} %
\DeclareMathOperator{\Ker}{Ker} %
\newcommand{\rest}[2]{#1\raisebox{-0.3ex}{\mbox{$\mid_{#2}$}}}

\begin{document}

\title[Cocycles and Ma\~{n}e sequences]{Cocycles and Ma\~{n}e sequences with an application to ideal fluids}
\author{Shvydkoy, Roman}
\address{University of Illinois at Chicago, Department of Mathematics
(M/C 249), Chicago, IL 60607} %
\email{shvydkoy@math.uic.edu}
\thanks{I thank Yuri Latushkin and Susan Friedlander for stimulating discussions.}%
\date{\today}%

\begin{abstract}
Exponential dichotomy of a strongly continuous cocycle $\bFi$ is
proved to be equivalent to existence of a Ma\~{n}e sequence either
for $\bFi$ or for its adjoint. As a consequence we extend some of
the classical results to general Banach bundles. The dynamical
spectrum of a product of two cocycles, one of which is scalar, is
investigated and applied to describe the essential spectrum of the
Euler equation in an arbitrary spacial dimension.
\end{abstract}

\maketitle

\section{Introduction}
Appearance of continuous cocycles on attractors of dissipative PDE,
in particular the Navier-Stokes equation, has spurred development of
the infinite-dimensional analogues of the classical results of
Mather \cite{Mather68}, Oseledets \cite{Os}, Sacker and Sell
\cite{SSSpec}, and others (see the historical account in \cite{CL}).
New implementations of the theory in ideal fluid dynamics raised new
questions that, as far as we know, have not been explicitly answered
before. One of them includes precise formulation of the relationship
between exponential dichotomy and existence of Ma\~{n}e sequences or
points. The purpose of this present note is partly to clarify known
results in this direction by proving them in the most general
settings, and partly to justify the cocycle related results claimed
in \cite{Shv2005}. The results have a direct physical interpretation
in terms of shortwave instabilities of an ideal fluid and apply to
describe the essential spectrum of non-dissipative advective
equations. One particular example we will consider in the next
section is the Euler equation on the torus. More examples and
further discussion can be found in these recent papers
\cite{FriedShv2005,FS2004,LatV2003,Shv2005}.

\section{Statements of the results}\label{S:state}

Let $\Th$ be a locally compact Hausdorff space countable at
infinity, and let $X$ be a Banach space. Suppose $\f = \{\f_t\}_{t
\in \R}$ is a continuous flow on $\Th$. A strongly continuous
exponentially bounded cocycle $\bFi$  over the flow $\f$ acting on
the trivial bundle $\Th \times X$ is a family of bounded linear
operators $\{ \bFi_t(\th)\}_{t \geq 0,\, \th \in \Th} \ss \cL(X)$
strongly continuous in $t,\th$, satisfying
$$
\bFi_0(\th) = \bI, \quad \bFi_t(\f_s(\th)) \bFi_s(\th) =
\bFi_{t+s}(\th),
$$
for all $\th \in \Th$, $t,s \geq 0$, and such that $ \sup_{0\leq t
\leq 1,\, \th \in \Th} \|\bFi_t(\th)\| < \infty$. Generic cocycles
appear as fundamental matrix solutions of systems of linear ODEs
with variable coefficients. So, the solution $f(t)$ of the Cauchy
problem
\begin{align}
 f_t & = \ba(\f_t(\th)) f, \label{ode1}\\
 f(0) & = f_0 \label{ode2}
\end{align}
is given by $f(t) = \bFi_t(\th) f_0$, where $\bFi$ is a cocycle over
$\f$.

Following Chow and Leiva \cite{ChowL96} we say that $\bFi$ has
\emph{exponential dichotomy} if there exists a continuous
projector-valued mapping $\bP(\th): X \ra X$ such that for some
$\e>0$ and $M>0$ one has
\begin{enumerate}
    \item $\bFi_t(\th) \bP(\th) = \bP(\f_t(\th)) \bFi_t(\th)$ ;
    \item $\sup_{\th \in \Th} \| \rest{\bFi_t(\th)}{\im \bP(\th)}\|
    \leq M e^{-\e t}$ ;
    \item the restriction $\rest{\bFi_t(\th)}{\Ker \bP(\th)} : \Ker \bP(\th) \ra  \Ker
    \bP(\f_t(\th))$ is invertible, and
    $$
    \|\bFi_t(\th) x\| \geq M^{-1} e^{\e t} \|x\|,
    $$
    holds for all $x \in \Ker \bP(\th)$, $t \geq 0$, and $\th \in \Th$.
\end{enumerate}

To every cocycle $\bFi$ we associate an evolution semigroup $\bE$ on
the space of $X$-valued continuous functions vanishing at infinity,
$C_0(\Th;X)$, acting by the rule
\begin{equation}\label{Mather}
    \bE_t f(\th) = \bFi_t(\f_{-t}(\th)) f(\f_{-t}(\th)).
\end{equation}
The following Dichotomy Theorem of Mather \cite{Mather68}, proved in
the general settings by Rau \cite{Rau}, and Latushkin and Schnaubelt
\cite{LatSchn99}, relates exponential dichotomy to the semigroup
$\bE$.
\begin{theorem}\label{T:dich}
The cocycle $\bFi$ has exponential dichotomy if and only if the
semigroup $\bE$ is hyperbolic on  $C_0(\Th;X)$, i.e. $\s(\bE_t) \cap
\T = \varnothing$, $t>0$.
\end{theorem}

Restatement of exponential dichotomy in terms of local growth
characteristics of the cocycle is our goal in this section. The
well-known lemma due to Ma\~{n}e says that in the case when $\Th$ is
compact, $\dim X<\infty$ and $\bFi$ is invertible, there exists a
point $\th_0\in\Th$ and vector $x_0\in X$ such that $\sup_{t\in \R}
\|\bFi_t(\th_0)x_0\| <\infty$, provided $1$ belongs to the
approximate point spectrum of $\bE_1$. Thus, by virtue of
\thm{T:dich}, if $\bFi$ has exponential dichotomy, then no such
(Ma\~{n}e) point and vector exist.

The analogue of Ma\~{n}e's lemma in the general settings was proved
in \cite{LatSchn99}, where points had to be replaced by so-called
Ma\~{n}e sequences.

\begin{definition}\label{D:mane} A sequence of pairs $\{(\th_n,x_n)\}_{n=1}^\infty$, where $\th_n \in \Th$
and $x_n \in X$, is called a \emph{Ma\~{n}e sequence} of the cocycle
$\bFi$ if $\{x_n\}_{n=1}^\infty$ is bounded and there are constants
$C>0$ and $c>0$ such that for all $n \in \N$
\begin{subequations}\label{e:mane1}
\begin{align}
\|\bFi_n(\th_n)x_n\| &>c, \label{e:mane1a}\\
\|\bFi_{k}(\th_n)x_n\| &<C, \text{ for all } 0\leq k\leq 2n.
\label{e:mane1b}
\end{align}
\end{subequations}
\end{definition}

In order to completely characterize the dichotomy in terms of
Ma\~{n}e sequences, one is lead to consider the adjoint operator
$\bE^*_1$ defined on the space of regular $X^*$-valued measures of
bounded variation, since if $\s_{ap}(\bE_1) \cap \T = \varnothing$,
then $\s_{p}(\bE^*_1) \cap \T \neq \varnothing$.

\begin{theorem}\label{T:char}
The following conditions are equivalent:
\begin{enumerate}
    \item[(i)] $\bFi$ is not exponentially dichotomic;
    \item[(ii)] There is a \Ms\ either for the cocycle $\bFi$ or for its
    adjoint $\bPsi$.
\end{enumerate}
\end{theorem}

We recall that the adjoint cocycle $\bPsi = \bFi^*$ is the cocycle
over the inverse flow $\{\f_{-t}\}_{t \in \R}$ acting on $\Th \times
X^*$ by the rule
$$
\bPsi_t(\th) = \bFi^*_t(\f_{-t}(\th)).
$$
It is the cocycle that generates the adjoint evolution semigroup
$\bE^*_t$, and it inherits the continuity and boundedness properties
from the original cocycle.

We note that in case $X$ is a Hilbert space, \thm{T:char} can be
deduced from the analogue of the Dichotomy Theorem \ref{T:dich} on
the space $L^2(\Th, \fm, X)$ over an appropriately chosen
$\f$-invariant measure $\fm$ (see \cite{Anton84,CS80,LatSt91}). In
this case one takes advantage of the apparent reflexivity of the
space and the $C^*$-algebra technique developed in
\cite{Antonevich}. As a corollary of \thm{T:char} and its proof we
will obtain the full analogue of the Dichotomy Theorem on $L^2(\Th,
\fm, X)$ for any, reflexive or not, Banach space $X$.

We now show an example of how \thm{T:char} applies to spectral
problems of fluid dynamics.

We consider the linearized Euler equation on the torus $\torus$:
\begin{align}
v_t &= - \der{u_0}{v} - \der{v}{u_0} - \n p, \label{EE}\\
\n \cdot v &= 0,\label{EE2}
\end{align}
where $u_0 \in [C^\infty(\torus)]^n$ is a given equilibrium solution
to the nonlinear equation. It can be shown that
\eqref{EE}--\eqref{EE2} generates a $C_0$-semigroup $\bG_t$ on the
space $L^2$ of divergence-free fields, and in fact, on any energy
Sobolev space $H^m$. In contrast to the point spectrum, the
essential spectrum of $\bG_t$ is related to so-called shortwave
instabilities of the flow $u_0$. Those are instabilities created by
localized highly oscillating disturbances of the form
$$
v_\d(x) = b_0(x) e^{i \xi_0 \cdot x /\d}, \quad \d \ll 1.
$$
Propagation of such disturbances along the corresponding streamline
of the flow $u_0$ can be described by the WKB-type asymptotic
formula
\begin{equation}\label{ans}
v(x,t) = b(x,t) e^{i S(x,t)/\d} + O(\d).
\end{equation}
In this formula the amplitude $b$ and frequency $\xi = \n S$ are
governed by evolution laws which can be obtained by direct
substitution of the ansatz \eqref{ans} into the linearized Euler
equation \eqref{EE}. In the Lagrangian coordinates associated with
the flow $u_0$, those laws become free of partial differentiation,
which allows one to view them as a finite-dimensional dynamical
system of the form \eqref{ode1}. Specifically, in this case $\Th =
\torus \times \rn$, $\f_t$ is the flow on $\Th$ generated by the
bicharacteristic system of equations describing evolution of the
material particle $x$ and frequency $\xi$:
\begin{align}
x_t &= u_0(x), \label{x} \\
\xi_t& = -\p u_0^{\top} \xi, \label{xi}
\end{align}
and the amplitude equation for $b(t)$ is given by
\begin{equation}\label{b}
b_t = \p u_0(x) b + \langle \p u_0(x) b ,\xi\rangle \xi
|\xi|^{-2},
\end{equation}
subject to incompressibility condition $b \perp \xi$ (see
\cite{Fincoll,Shv2005} for details).

Let $\bB$ stand for the cocycle generated by the amplitude equation
\eqref{b}, and let $\chi_t$ denote the integral flow of $u_0$, i.e.
the solution of \eqref{x}. It can be shown that in terms of $\bB$
and $\chi$ the asymptotic formula \eqref{ans} takes the form
\begin{equation}\label{asym}
v(x,t) = \bG_t v_\d (x) = \bB_t(\chi_{-t}(x),
\xi_0)v_0(\chi_{-t}(x)) + O(\d),
\end{equation}
as $\d \ra 0$ (see \cite{Shv2005,V96}). Thus, if the cocycle $\bB$
has growing solutions, then the semigroup $\bG_t$  and hence the
flow $u_0$ is linearly unstable to shortwave perturbations.

Suppose now that $\bB$ is not dichotomic. Then in view of
\thm{T:char} either $\bB$ or $\bB^*$ has a \Ms. Since $\bG_t$
corresponds to $\bB^*$ through a formula similar to \eqref{asym},
and $\s(\bG_t) = \s(\bG_t^*)$, we can assume for definiteness that
$\bB$ has a \Ms\ $\{(x_n,\xi_n),b_n\}$. We consider a vector field
$b_n(x)$ localized near $x_n$ and aligned with $b_n$ up to a term of
order $O(\d)$ so that $v_{\d,n} =b_n(x) e^{i \xi_n \cdot x /\d}$ is
divergence-free. Choosing $\d = \d_n$ small enough we obtain
\begin{equation}\label{Gk}
\bG_k v_{\d,n} =\bB_k(\chi_{-k}(x), \xi_n)b_n(\chi_{-k}(x))e^{i
\xi_n \cdot \chi_{-k}(x) /\d} + O(\d),
\end{equation}
for all $ 0 \leq k \leq 2n$. Thus, denoting $z_n = v_{\d(n),n}$ we
fulfill the sufficient condition for hyperbolicity stated in
\lem{L:hypo} below, which implies that $1 \in |\s(\bG_t) |$. Given
the fact that the constructed sequence $z_n$ is weakly-null, we can
even conclude that $1 \in |\s_{ess}(\bG_t) |$, the essential
spectrum in the Browder sense.

\begin{definition}
Let us recall that the \emph{dynamical spectrum} of a cocycle $\bFi$
is the set of all points $\l \in \R$ such that  $\{ e^{-\l t} \bFi_t
\}$ has no exponential dichotomy. We denote this set by
$\Sigma_\bFi$.
\end{definition}
Generally, for a cocycle $\bFi$ with compact fiber-maps, its
dynamical spectrum is the union of disjoint segments, which may tend
to $-\infty$ or be infinite on the left (see \cite{CL,Mg2,SSSpec}).
Moreover, the number of segments is limited to the spacial dimension
of $X$ if the latter is finite.

After rescaling, \thm{T:char} states that $\l \in \Sigma_\bFi$ if
and only if either $e^{-\l t} \bFi_t$ or its adjoint has a \Ms.
Thus, going back to our example we obtain the following inclusion
\begin{gather}
\exp\{ t \Sigma_\bB \} \ss |\s_{ess}(\bG_t) |, \\
\intertext{while, on the other hand, as shown in
\cite{Shv2005,V96},}
|\s_{ess}(\bG_t) | \ss \exp\{ t [\min
\Sigma_\bB,\max \Sigma_\bB]\}.
\end{gather}
In view of the above discussion the physical meaning of a \Ms\
becomes more transparent in the context of fluid dynamics: it shows
exactly what particle in what frequency has to be excited to
destabilize the flow. The dynamical spectrum $\Sigma_\bB$, in turn,
provides the range of all possible rates at which the excitations
grow exponentially.

On the Sobolev space $H^m$ of divergence-free fields, the norm of
$\bG_k v_{\d,n}$  behaves like $\| \p \chi_k^{-\top}(x_n) \xi_n\|^m
\| \bB_k(x_n,\xi_n) b_n \|$, as $\d \ra 0$. So, in this case one is
naturally lead to consider the augmented cocycle
$$
\bB\bX^m_t (x,\xi) = \| \p \chi_t^{-\top}(x) \xi \|^m  \bB_t(x,\xi).
$$
Via a similar reasoning as above we can obtain the following
inclusions
\begin{equation}\label{sobincl}
\exp\{ t \Sigma_{\bB\bX^m} \} \ss |\s_{ess}(\bG_t) | \ss \exp\{ t
[\min \Sigma_{\bB\bX^m},\max \Sigma_{\bB\bX^m}]\}.
\end{equation}

The influence of the scalar cocycle $\bX_t^m =\| \p \chi^{-\top}(x)
\xi\|^m$ on the whole spectrum of $\bB\bX^m$ is growing with $m$
provided $\bX^m$ itself has a non-trivial spectrum, or equivalently,
$u_0$ has exponential stretching of trajectories. Since
$\Sigma_\bX^m$ is one connected segment expanding as $m \ra \infty$,
it will fill all possible gaps in $\Sigma_{\bB\bX^m}$ for $m$ large
enough. Whenever this happens we obtain the identity
\begin{equation}\label{sobident}
\exp\{ t \Sigma_{\bB\bX^m} \} = |\s_{ess}(\bG_t) |.
\end{equation}
According to \thm{T:char}, to every point of the set
$\Sigma_{\bB\bX^m}$ there corresponds a \Ms. The fact that this set
gets larger with $m$ and eventually becomes connected implies that
there is an increasing number of \Ms\ needed to serve points of
$\Sigma_{\bB\bX^m}$. Physically, this means that in a finer norm,
such as the norm of $H^m$, fluid has more spots sensitive to
shortwave perturbations than it does in the basic energy norm.
Although the above statements apply in any spacial dimension, in the
more tractable case of $n=2$ a much stronger result was obtained by
Koch \cite{Koch2002}. It shows that any non-isochronic stationary
flow in a flat domain is nonlinearly instable in the H\"{o}lder
classes $C^{1,\a}$.

Motivated by the example of the Euler equation, in Section
\ref{S:scalar} we pose the general question of how the spectrum of a
cocycle $\bFi$ changes under multiplication by another scalar
cocycle $\bC$. We will show that it is contained in the arithmetic
sum of $\Sigma_\bFi$ and $\Sigma_\bC$, and we give a sufficient
condition for $\Sigma_{\bC\bFi}$ to be connected. This condition
applied to the Euler equation will yield a lower bound on $m$ for
which \eqref{sobident} holds. This will completely justify the
result claimed in \cite{Shv2005}.

Finally, we remark that all our arguments are local, and as such can
be generalized to an arbitrary continuous Banach bundle.

\section{Characterization of exponential dichotomy}\label{S:char}

In this section we present the proof of \thm{T:char} and use it to
show the analogue of the Dichotomy Theorem on $L^p$ spaces.

The proof relies on the following lemma, which we state slightly
more generally than it is needed at the moment. However, it will be
used later to its full extent.

\begin{lemma}\label{L:hypo} Let $Z$ be a Banach space and $\bT \in
\cL(Z)$. Suppose there is a bounded sequence of vectors $\{z_k\}_{k
= 1}^\infty$ and a subsequence of natural numbers
$\{n_k\}_{k=1}^\infty \ss \N$ such that
\begin{itemize}
    \item[(a)] $\lim_{k \ra \infty} n_k^{-1} \log \| \bT^{n_k} z_k \|
    \geq \l_1$;
    \item[(b)] $\lim_{k \ra \infty} n_k^{-1} \log \| \bT^{2n_k} z_k \|
    \leq \l_1+ \l_2$.
\end{itemize}
Then the following statements are true
\begin{itemize}
    \item[(i)] If $\l_1 \leq \l_2$, then $[\l_1,\l_2] \cap \log |
    \s(\bT)| \neq \varnothing$;
    \item[(ii)] If $\l_2 \leq \l_1$, then $[\l_2,\l_1] \ss \log |
    \s(\bT)|$.
\end{itemize}
\end{lemma}
\begin{proof}
To prove (i) let us assume, on the contrary, that $[\l_1,\l_2] \cap
\log|\s(\bT)| =\varnothing$. Then there is $\e >0$ such that
$$
[\l_1-\e,\l_2+\e] \cap \log|\s(\bT)| =\varnothing.
$$
Let $Z_s$ and $Z_u$ denote the spectral subspaces corresponding to
the parts of the spectrum below $\l_1 - \e$ and above $\l_2 +\e$,
respectively. For $n$ large enough we have
\begin{align}
\| \rest{\bT^n}{Z_s} \| & < e^{n(\l_1 - \e)}; \label{hyp1}\\
\| \bT^n z \| & \geq  e^{n(\l_2+ \e)} \|z\|, \quad z \in Z_u.
\label{hyp2}
\end{align}
Let $z_k = z_k^s + z_k^u$. Then, by (b) and \eqref{hyp2},
\begin{multline*}
\l_1 + \l_2 \geq \lim_{k \ra \infty} n_k^{-1} \log \| \bT^{2n_k}
z_k^u \| \geq \l_2 + \e + \lim_{k \ra \infty} n_k^{-1} \log \|
\bT^{n_k} z_k^u \|.
\end{multline*}
So,
$$
\lim_{k \ra \infty} n_k^{-1} \log \| \bT^{n_k} z_k^u \| \leq \l_1 -
\e.
$$
In combination with \eqref{hyp1} this gives
$$
\lim_{k \ra \infty} n_k^{-1} \log \| \bT^{n_k} z_k \| \leq \l_1 -\e,
$$
which contradicts condition (b).

To prove (ii) let us assume that $\l_2 < \l_1$ and fix any $\l \in
[\l_2,\l_1]$. Let us denote $\d = \l_1 - \l \geq 0$. We consider a
new bounded sequence
$$
w_k = e^{-n_k \d} z_k.
$$
For this sequence the following conditions are verified:
\begin{align*}
\lim_{k \ra \infty} n_k^{-1} \log \| \bT^{n_k} w_k \| & \geq \l, \\
\lim_{k \ra \infty} n_k^{-1} \log \| \bT^{2n_k} w_k \| & \leq \l_2
+\l \leq 2\l.
\end{align*}
Applying (i) with $\l_1 = \l_2 = \l$ we obtain $\l \in
\log|\s(\bT)|$.
\end{proof}

\begin{proof}[Proof of \thm{T:char}]

Let us assume (ii). Suppose $\{(\th_n,x_n)\}_{n=1}^\infty$ is a \Ms\
for the cocycle $\bFi$. For each $n$ let us find an open
neighborhood of $\th_n$, denoted $U_n$, such that for all $\th \in
U_n$,
\begin{subequations}\label{subineq}
\begin{align}
\|\F_n(\th_n) - \F_n(\th)\| &<c/2 ,\\
\|\F_{2n}(\th)\| &< 2C.
\end{align}
\end{subequations}

Let $\phi_n \in C_0(\Th)$ be a scalar function of unit norm
supported on $U_n$ such that $\phi_n(\th_n) = 1$. Then, by
\eqref{subineq} and \eqref{e:mane1}, we obtain
\begin{align*}
\|\bE_n(\phi_n(\cdot) x_n )\| & > c/2 ,\\
\|\bE_{2n}(\phi_n(\cdot) x_n)\| &< 2C.
\end{align*}
So, Lemma \ref{L:hypo} applies with $\l_1=\l_2 = 0$ to show that
$\bE$ is not hyperbolic, and hence, by virtue of \thm{T:dich},
$\bFi$ is not exponentially dichotomic.

If there is a \Ms\ for the adjoint cocycle $\bPsi$, then by the
previous argument applied to $\bE^*$ on $C_0(\Th; X^*)$, we find
that $\bPsi$ is not exponentially dichotomic. Hence, $\bFi$ is not
dichotomic either, as seen directly from the definition.

To show the converse implication, let us assume (i). By Theorem
\ref{T:dich} one has $1 \in | \s(\bE_1)|$. There are two
possibilities that follow from this -- either there is an
approximate eigenvalue or there is a point of the residual spectrum
on the unit circle.

In the first case there is a normalized sequence of functions $f_n
\in C_0(\Th;X)$ such that
\begin{equation}\label{appr}
\|\bE_k f_n - e^{i\a k} f_n \| \leq \frac{1}{2},
\end{equation}
for some $\a \in \R$ and all $k = 1,\ldots,2n$. Let us choose points
$\th'_n \in \Th$ so that $\|f_n(\th'_n)\| = 1$. By \eqref{appr}, we
have
\begin{gather*}
\| \F_n(\f_{-n}(\th'_n)) f_n(\f_{-n}(\th'_n))\| \geq 1/2,\\
\intertext{and} %
\| \F_{k}(\cdot)f_n(\cdot) \| \leq 2, \quad 1 \leq k \leq 2n.
\end{gather*}
Choosing $\th_n =\f_{-n}(\th'_n)$ and $x_n =f_n(\f_{-n}(\th'_n))$ we
fulfill the conditions of \defin{D:mane}.

In the second case, let $e^{i \a}$ be a point of the residual
spectrum of $\bE_1$. Hence, there exists $\nu \in \M(\Th, X^*)$,
with $\|\nu\| = 1$, a regular Borel $X^*$-valued measure of bounded
variation, such that
\begin{equation}\label{eigenmeasure}
\bE_n^* \nu = e^{i \a n} \nu, \quad n \in \N.
\end{equation}

Recall that the norm in $\M(\Th;X^*)$ is given by the total
variation
\begin{equation}\label{totvar}
    \|\nu\| = \sup \left\{ \sum_{i=1}^N \|\nu(A_i)\|:
    \bigcup_{i=1}^N A_i  = \Th,\ A_i \cap A_j = \emptyset
    \right\}.
\end{equation}
We also consider the semivariation of a set $A \ss \Th$ defined by
$$
|\nu|(A) = \sup\{ |x^{**}\nu|(A): x^{**} \in X^{**}\},
$$
and we recall the following inequality \cite[p.4]{Diestel-Uhl}:
\begin{equation}\label{sv}
    |\nu|(A) \leq 4 \sup\{\|\nu(B)\|: B \ss A\}.
\end{equation}

Going back to our proof, let us fix $n \in \N$. By the continuity of
$\bPsi$ and $\f$, using the topological assumption on $\Th$, we can
find a partitioning of $\Th$ into Borel sets $\{A_j\}_{j \in J}$
such that for every $j \in J$,
\begin{align}
   \| \bPsi_{n}(\f_{-n}(\th')) - \bPsi_{n}(\f_{-n}(\th'')) \| & < c_0,
   \label{close1} \\
   \| \bPsi_{2n}(\f_{-n}(\th')) - \bPsi_{2n}(\f_{-n}(\th'')) \| & <
   c_0, \label{close2}
\end{align}
holds for all $\th', \th'' \in A_j$, and where the constant $c_0>0$
is to be specified later.

By \eqref{sv}, for every $j \in J$, there is a set $B_j \ss A_j$
such that
\begin{equation}\label{Bj}
4 \| \nu(B_j)\| > |\nu|(A_j).
\end{equation}
Let us fix arbitrary tag points $\th_j \in B_j$. According to
\eqref{eigenmeasure} and \eqref{close1} -- \eqref{close2}, we have
\begin{align}
\bPsi_{n}(\f_{-n}(\th_j)) \frac{\nu(\f_{-n}(B_j))}{|\nu|(A_j)} & =
\frac{\nu(B_j)}{|\nu|(A_j)} + v_n^j, \label{3}\\
\bPsi_{2n}(\f_{-n}(\th_j)) \frac{\nu(\f_{-n}(B_j))}{|\nu|(A_j)} & =
\frac{\nu(\f_n(B_j))}{|\nu|(A_j)} + u_{n}^j, \label{4}\\
 \intertext{where}
 \|v_{n}^j\|, \|u_n^j\| &< c_0
 \frac{\|\nu(\f_{-n}(B_j))\|}{|\nu|(A_j)}. \label{5}
\end{align}

Let us denote $\eta = |\nu|(\Th)$. We claim that there exists $j =
j(n) \in J$ such that
\begin{align}
\|\nu(\f_{-n}(B_{j(n)}))\| &\leq \frac{4}{\eta} |\nu| (A_{j(n)}),\label{6}\\
\|\nu(\f_{n}(B_{j(n)}))\| &\leq \frac{4}{\eta} |\nu|
(A_{j(n)}).\label{7}
\end{align}

Indeed, suppose there is no such $j(n)$. Then for each $j \in J$
either \eqref{6} or \eqref{7} fails. So, by the subadditivity of
semivariation, we obtain
\begin{multline*}
\eta = |\nu|(\Th) \leq \sum_{j \in J} |\nu|(A_j) \leq
\frac{\eta}{4}\sum_{j\in J}\|\nu(\f_{-n}(B_{j}))\|
+\\+\frac{\eta}{4}\sum_{j\in J}\|\nu(\f_{n}(B_{j}))\| \leq 2
\frac{\eta}{4} \|\nu\| = \frac{\eta}{2},
\end{multline*}
a contradiction.

Let us put $\th_n = \f_{-n}(\th_{j(n)})$ and $x^*_n =
\frac{\nu(\f_{-n}(B_{j(n)}))}{|\nu|(A_{j(n)})}$. In view of
\eqref{6}, $\{x_n^*\}$ is a bounded sequence. Also, by \eqref{6} and
\eqref{5}, we have
$$
\|v_{n}^{j(n)}\|, \|u_{n}^{j(n)}\| < \frac{4 c_0}{\eta}.
$$
So, by  \eqref{Bj}, \eqref{3}, \eqref{4}, \eqref{7}, and \eqref{sv},
\begin{align}
\|\bPsi_n(\th_n) x^*_n\| \geq \frac{1}{4} - \frac{4c_0}{\eta},\\
\|\bPsi_{2n}(\th_n)x_n^*\| \leq \frac{4}{\eta} + \frac{4c_0}{\eta}.
\end{align}
It suffices to take $c_0 = \eta /32$.

\end{proof}

In the compact case existence of a Ma\~{n}e sequence is equivalent
to existence of a Ma\~{n}e point (see, for example, \cite{CL}). So,
in this case \thm{T:char} can be restated as follows.

\begin{corollary}
Suppose $\dim X < \infty$ and $\Th$ is compact. Then $\bFi$ is
exponentially dichotomic if and only if either $\bFi$ or $\bPsi$ has
a Ma\~{n}e point.
\end{corollary}

Another fact that follows directly from \lem{L:hypo} is that any
Lyapunov index of the cocycle $\bFi$ belongs to the dynamical
spectrum $\Sigma_\bFi$ (see also Johnson, Palmer and Sell
\cite{JPS87}).

Indeed, suppose
$$
\l = \lim_{k \ra \infty} n_k^{-1} \log \| \bFi_{n_k} (\th) x \|,
$$
for some $\th \in \Th$ and $x \in X$. Then by the same construction
as in the proof of \thm{T:char} we find functions $f_k$ such that
\begin{align*}
\l &= \lim_{k \ra \infty} n_k^{-1} \log \| \bE_{n_k} f_k \|, \\
2\l &= \lim_{k \ra \infty} n_k^{-1} \log \| \bE_{2n_k} f_k \|.
\end{align*}
Applying \lem{L:hypo} with $\l_1 = \l_2 = \l$ we obtain $\l \in \log
|\s(\bE_t)|$.

As another consequence of \thm{T:char} we prove the analogue of the
Dichotomy Theorem \ref{T:dich} for $L^p$-spaces.

Let $\fm$ be a Borel $\f$-quasi-invariant measure on $\Th$. We
define $\bE$ on $L^p(\Th, \fm, X)$, $1\leq p<\infty$, by the rule
\begin{equation}\label{MatherLp}
    \bE_t f(\th) = \left(  \frac{ d(\fm \circ \f_{-t})}{d\fm} \right)^{1/p}\bFi_t(\f_{-t}(\th))
    f(\f_{-t}(\th)),
\end{equation}
where the expression under the root is the Radon-Nikodim derivative
(we refer to \cite{CL} for a detailed discussion).

\begin{theorem}\label{T:charlp} Let $\bE$ be defined by
\eqref{MatherLp} on the space $L^p(\Th, \fm, X)$, with $1\leq
p<\infty$, where $\fm$ is a Borel $\f$-quasi-invariant measure such
that $\fm(U) > 0$ for every open set $U$. Then $\bFi$ has
exponential dichotomy if and only if\, $\bE$ is hyperbolic.
\end{theorem}
\begin{proof}
Suppose that $\bFi$ has exponential dichotomy, then the spaces
\begin{align}
Z_s & = \{ f \in L^p(\Th, \fm, X): f(\th)\in \im \bP(\th) \}, \\
Z_u & = \{ f \in L^p(\Th, \fm, X): f(\th)\in \Ker \bP(\th) \}
\end{align}
define, respectively, exponentially stable and unstable subspaces
for $\bE$ such that $L^p(\Th, \fm, X) = Z_s \oplus Z_u$. Hence,
$\bE$ is hyperbolic.

Suppose $\bFi$ has no exponential dichotomy. Let us assume that
$\bFi$ has a \Ms. Then the same construction as in the proof of
\thm{T:char}, with localized scalar functions $\phi_n \in
L^p(\Th,\fm)$, $\|\phi_n\|_p = 1$, shows that $\bE$ is not
hyperbolic.

If the adjoint cocycle $\bPsi$ has a Ma\~{n}e sequence, then we
regard the corresponding functions $\f_n(\th) x_n^*$ as elements of
$L^q_{w^*}(\Th,\fm,X^*)$, the space of weak$^*$-measurable
$q$-integrable functions with values in $X^*$. This space is the
dual of $L^p(\Th,\fm,X)$, provided $p^{-1} + q^{-1} = 1$ (see
\cite{VVbook}).

From \lem{L:hypo} we conclude that the operator $\bE^*_1$ is not
hyperbolic over $L^q_{w^*}(\Th,\fm,X^*)$. Hence, $\bE_1$ is not
hyperbolic over $L^p(\Th,\fm,X)$.
\end{proof}

\section{Scalar multiple of a cocycle}\label{S:scalar}

Let $\Th,\f,X$, and $\bFi$ be as before, and let $\bC =
\{\bC_t(\th)\}_{t \geq 0,\, \th \in \Th}$ be a scalar cocycle over
the same flow $\f$ acting on $\Th \times \C$. Then the product
$\bC\bFi =\{\bC_t(\th)\bFi_t(\th)\}_{t \geq 0,\, \th \in \Th}$
defines another cocycle on $\Th \times X$. An example of how
products of this type arise in the equations of fluid dynamics was
presented in Section \ref{S:state}.

\begin{lemma}\label{L:incl}
One has the following inclusion
\begin{equation}\label{E:incl}
    \Sigma_{\bC\bFi} \ss \Sigma_\bC + \Sigma_{\bFi}.
\end{equation}
\end{lemma}
\begin{proof}
Let $\rho \in \Sigma_{\bC\bFi}$. Then by \thm{T:char} there exists a
\Ms, say, for $e^{-\rho t} \bC_t\bFi_t$ (the case of adjoint cocycle
is treated similarly). Let $\{\th_n,x_n\}_{n=1}^\infty$ be that
sequence. Then we have
\begin{align}
|\bC_n(\th_n)| \| \bFi_n(\th_n) x_n \| & > c e^{\rho t}, \label{l1} \\
|\bC_{2n}(\th_n)| \| \bFi_{2n}(\th_n) x_n \| & < C e^{2 \rho t},
\label{l2}
\end{align}
for all $n \in \N$. Let us extract a subsequence
$\{n_k\}_{k=1}^\infty$ such that the limits
\begin{align}
\lim_{k \ra \infty} n_k^{-1} \log |\bC_{n_k}(\th_{n_k}) | & = \l_1,\label{l3} \\
\lim_{k \ra \infty} n_k^{-1} \log |\bC_{2n_k}(\th_{n_k}) | & = \l_1+\l_2, \label{l4}\\
\lim_{k \ra \infty} n_k^{-1} \log \|\bFi_{n_k}(\th_{n_k})x_{n_k} \| & = \mu_1, \label{l5}\\
\lim_{k \ra \infty} n_k^{-1} \log \|\bFi_{2n_k}(\th_{n_k})x_{n_k} \|
& = \mu_1+\mu_2, \label{l6}
\end{align}
exist. By \eqref{l1} and \eqref{l2}, we have
\begin{align}
\l_1 + \mu_1 \geq \rho, \label{l7}\\
\l_1 +\l_2 + \mu_1 +\mu_2 \leq 2 \rho. \label{l8}
\end{align}

Let us consider two cases: $\l_1 \leq \l_2$ and $\l_1 > \l_2$. If
$\l_1 \leq \l_2$, then by \lem{L:hypo}, there is $\l \in [\l_1,\l_2]
\cap \Sigma_\bC$. From \eqref{l7} and \eqref{l8}, we have
$\l_2+\mu_2 \leq \rho$. So, $\mu_2 \leq \rho - \l_2 \leq \rho - \l_1
\leq \mu_1$. In this case \lem{L:hypo} implies that $[\mu_2,\mu_1]
\ss \Sigma_{\bFi}$. We choose $\mu = \rho - \l \in [\mu_2,\mu_1]$ to
satisfy $\rho = \l + \mu$.

If $\l_2 < \l_1$, then $[\l_2, \l_1] \ss \Sigma_\bC$. From the above
we still have $\mu_2 \leq \rho - \l_2$ and $\mu_1 \geq \rho - \l_1$.
If $\mu_1 \leq \mu_2$, then we find a point $\mu \in [\mu_1,\mu_2]
\cap \Sigma_{\bFi}$, and choose $\l = \rho - \mu \in [\l_2, \l_1]$.
If $\mu_1 > \mu_2$, then $[\mu_2,\mu_1] \ss \Sigma_{\bFi}$ and
$[\mu_2,\mu_1] \cap [\rho-\l_1,\rho-\l_2] \neq \varnothing$.
Choosing $\mu \in [\mu_2,\mu_1] \cap [\rho-\l_1,\rho-\l_2]$ we get
$\l = \rho - \mu \in [\l_2, \l_1]$. This finishes the argument.

\end{proof}

\def \fmax {\mu_{\mathrm{max}}^\bFi}
\def \fmin {\mu_{\mathrm{min}}^\bFi}
\def \cmax {\mu_{\mathrm{max}}^\bC}
\def \cmin {\mu_{\mathrm{min}}^\bC}
\def \cfmax {\mu_{\mathrm{max}}^{\bC\bFi}}
\def \cfmin {\mu_{\mathrm{min}}^{\bC\bFi}}

Now let us assume that both cocycles $\bFi$ and $\bC$ are invertible
so that their spectra are bounded from above and below. We denote
$\fmax = \max \Sigma_{\bFi}$ and $\fmin = \min \Sigma_\bFi$. Similar
notation will be used for other cocycles.

\begin{lemma} \label{L:res}
Suppose $\rho \in [\cfmin,\cfmax] \backslash \Sigma_{\bC\bFi}$. Then
the following inequalities hold:
\begin{equation}\label{E:res}
    \cmax + \fmin < \rho < \cmin + \fmax.
\end{equation}
\end{lemma}
\begin{proof} Let $\bP$, $\e>0$ and $M$ be as in the definition of
the dichotomy. Let us fix any $\l \in \Sigma_{\bC}$. Then by
\thm{T:char} there exists a \Ms\ $\{\th_n\}$ for $\bC$:
\begin{align}
|\bC_n(\th_n)| &> c e^{n \l}, \label{c1}\\
|\bC_{2n}(\th_n)| &< C e^{ 2n \l}. \label{c2}
\end{align}

Given that $\Ker \bP(\th) \neq \{0\}$ for every $\th \in \Th$, we
can find a unit vector $x_n \in \Ker \bP(\th_n)$ for every $n$. Then
in view of \eqref{c2} we have
$$
C e^{n (2\l + 2\fmax + \e)} \geq |\bC_{2n}(\th_n)| \|
\bFi_{2n}(\th_n) x_n\| \geq M^{-1} e^{2n(\rho + \e)}.
$$
Thus, $\rho \leq \l + \fmax-\e$ for all $\l \in \Sigma_\bC$. This
proves the right side of \eqref{E:res}. The left side is proved
similarly using \eqref{c1}.
\end{proof}

As an immediate consequence of \lem{L:res} we obtain the following
sufficient condition for $\Sigma_{\bC\bFi}$ to be connected.

\begin{theorem}\label{T:conn}
Suppose the cocycles $\bC$ and $\bB$ are invertible. The dynamical
spectrum $\Sigma_{\bC\bFi}$ has no gaps provided the diameter of
$\Sigma_{\bC}$ is greater than the diameter of $\Sigma_{\bFi}$, i.e.
\begin{equation}\label{E:conn}
    \cmax - \cmin \geq \fmax - \fmin.
\end{equation}
\end{theorem}

\def \lmax {\l_{\mathrm{max}}}
\def \lmin {\l_{\mathrm{min}}}

Going back to our example with the Euler equation, let us denote
$$
\lmax =\mu_{\mathrm{max}}^\bX, \quad \lmin =\mu_{\mathrm{min}}^\bX.
$$
Then $\Sigma_{\bX^m} = m [\lmin,\lmax]$. Assume that $\lmax >0$, and
hence by incompressibility, $\lmin <0$. In this case condition
\eqref{E:conn} turns into
\begin{equation}\label{Eulconn}
    |m| \geq \frac{\mu_{\mathrm{max}}^\bB -
    \mu_{\mathrm{min}}^\bB}{\lmax - \lmin}.
\end{equation}
So, if $|m|$ is large enough, then we have identity \eqref{sobident}
over the Sobolev space $H^m$. In fact, if the cocycle $\bB$ has
trivial dynamical spectrum, such as in the case of a parallel shear
flow $u_0$ or $n=2$ in the vorticity formulation, then
$\mu_{\mathrm{max}}^\bB =\mu_{\mathrm{min}}^\bB$, and
\eqref{sobident} holds for any $m \neq 0$.

We refer to \cite{SL2003b,Shv2005} for more details on the
description of the essential spectrum for the Euler and other
similar equations.


\def\cprime{$'$} \def\cprime{$'$} \def\cprime{$'$} \def\cprime{$'$}
  \def\cprime{$'$} \def\cprime{$'$} \def\cprime{$'$} \def\cprime{$'$}
  \def\cprime{$'$} \def\cprime{$'$} \def\cprime{$'$} \def\cprime{$'$}

\end{document}